\documentclass[11pt,a4paper]{amsart}

\usepackage{epsf}

\title{    
Non-connected toric Hilbert schemes
}
\author{     Francisco Santos}
\address{    Departamento de Matem\'aticas, Estad\'{\i}stica y Computaci\'on\\
             Universidad de Cantabria, E-39005, Santander, SPAIN.}

\email{santosf@unican.es
}

\thanks{This result was obtained in the fall of 2001,
while I was a visiting professor in the
Department of Mathematics, U.C. Davis, supported by U. C. Davis,
M.S.R.I. and the spanish government. I am also partially supported by grant
BFM2001--1153 of the Spanish  Direcci\'on General de Ense\~nanza 
Superior e Investigaci\'on Cient\'{\i}fica}

\date{April 2002, revised November 2003. 
Dedicated to Bernd Sturmfels on his 40th birthday.}

\newtheorem{theoremintro}{Theorem}
\newtheorem{theorem}{Theorem}[section]
\newtheorem{lemma}[theorem]{Lemma}
\newtheorem{corollary}[theorem]{Corollary}
\newtheorem{proposition}[theorem]{Proposition}

\theoremstyle{definition}
\newtheorem{defi}[theorem]{Definition}

\theoremstyle{remark}
\newtheorem{remark}[theorem]{Remark}

\newcommand{\noproof}{\qed}

\newcommand{\naturals}{{\mathbb N}}
\newcommand{\integers}{{\mathbb Z}}
\newcommand{\reals}{{\mathbb R}}
\newcommand{\field}{{k}}
\newcommand{\rationals}{{\mathbb Q}}

\newcommand{\conv}{{\hbox{\rm conv}}}
\newcommand{\Rad}{{\hbox{\rm Rad}}}

\newcommand{\bb}{{\mathbf b}}
\newcommand{\cc}{{\mathbf c}}
\newcommand{\e}{{\overline e}}
\renewcommand{\aa}{{\mathbf a}}

\newcommand{\A}{{\mathcal A}}

\newcommand{\K}{{\mathcal K}}

\begin{document}

\begin{abstract}
We construct small (50 and 26 points, respectively) point sets
in dimension 5 whose graphs of triangulations are  not connected. 
These examples improve
our construction in \textit{J. Amer. Math. Soc.} \textbf{13}:3 (2000),
611--637 not only in size, but also in that the associated toric Hilbert schemes 
are not connected either, a question left open in that article.
Additionally, the point sets can easily
be put into convex position, providing examples
of 5-dimensional polytopes with non-connected graph of triangulations.
\end{abstract}

\maketitle

\section*{\bf Introduction}
\addcontentsline{toc}{section}{Introduction}
\label{intro}

The graph of triangulations of a finite point set $A\subset\reals^d$ 
has as vertices all the triangulations of $A$ and as edges certain natural local
modifications of them, analogous to the
bistellar flips considered in combinatorial topology \cite{BjoLut}.
See the precise definition in Section
\ref{sec.prelim}. 
This graph is interesting in several contexts:

In Geometric Combinatorics, its study is a special case of the 
  Generalized Baues problem, which includes
  monotone paths on polytopes, zonotopal tilings, or the extension space of
  realizable oriented matroids
as other cases. See the survey \cite{Reiner-survey}, or
\cite{BiKaSt,Santos-refine}.

In Computational Geometry, triangulations are a standard tool and flips
  have been often proposed as a method to explore the space of all possible
  triangulations or search for optimal ones (e.g., the Delaunay one). See
  \cite{EdeSha,HuNoUr,Lawson} or the survey \cite{BerEpp}.

 In Algebraic Geometry, lattice polytopes and triangulations of them
  are closely related to toric 
  varieties~\cite{GKZbook,IteRoy,Sturmfels-book}.
 In particular, 
 the graph of triangulations of a point set $A$ is connected if and only
 if the Chow variety \cite{KaStZe}
 of $A$ is connected (see Corollary 4.9 in \cite{Santos-noflips}). 
 The same question for the  \emph{toric Hilbert scheme} 
 \cite{MacTho,PeeSti,StStTh} of the point set
 is not clear. Even though
 Sturmfels and Haiman \cite{HaiStu} have recently constructed a natural
 morphism from the toric Hilbert scheme to the toric Chow
 variety, this morphism is in general not surjective. In particular,
 it remained open until the present construction whether non-connected
toric Hilbert schemes existed. The solution we present here uses the result,
due to Maclagan and Thomas \cite{MacTho}, that at least the
 subgraph induced by unimodular triangulations of $A$ has a
faithful relation to the toric
 Hilbert scheme. See  Section \ref{sub.hilbert}.

More details on these interrelations 
can also be found in \cite{Santos-noflips}, where we
constructed the first example of a point set 
whose graph of triangulations is not
connected, in dimension 6 with 324 points.
Here we show much smaller examples. The essential new
idea is that we construct
triangulations with sub-structures which cannot be destroyed by flips, rather
than triangulations without flips.
This gives more flexibility to the construction, and it
also saves us several technicalities in the proofs.

The paper is organized as follows: in Section 
\ref{sec.subcomplexes}  we give the precise definition
of the graph of triangulations and the main ingredients needed to
prove that the graph is not connected in our examples, summarized in Theorem
\ref{thm.map} and its Corollary \ref{coro.map}. 
This section describes also
the connection between the graph of triangulations and the toric Hilbert
scheme.

The other two sections
describe our two point sets $A_{50}$ and $A_{26}$ and prove that their graphs
of triangulations and their toric Hilbert schemes
are not connected. These point sets have 50 and 26 points
respectively, both in dimension 5. The reason why we include both, and not
only the smaller one, is that $A_{50}$ is simpler to describe and 
it has stronger
properties which make it easier to prove what we need. More precisely,
\[
A_{50} = \left(A\cup\{O\}\right)\times \{0,1\}
\]
where $A$ and $O$ are, respectively, the set of vertices
and the centroid of the
regular 24-cell in $\reals^4$. The point 
set $A_{26}$ is related to the regular
cross-polytope, but not in such a direct way.

Our main results (Theorem  \ref{thm.24cell}, 
Theorem \ref{thm.26points} and Propositions
\ref{prop.50convex} and  \ref{prop.26convex})
can be summarized as:

\begin{theoremintro}
\label{thm.main}
The 5-dimensional point sets $A_{26}$ and $A_{50}$, with 26 and 50 points
respectively, have the following properties:
\begin{enumerate}
\item Their graphs of triangulations are not connected. 

\item Their toric Hilbert schemes are not connected (for $A_{26}$, the 
  scheme is considered with 
  respect to a non-homogeneous grading).

\item The graphs of triangulations remain non-connected under suitable
  perturbations of the point sets into convex position. In particular, 
  there is a
  5-dimensional polytope with 26 vertices whose graph of triangulations
is not
  connected. 
\end{enumerate}
\end{theoremintro}

Concerning part (3),
a lifting construction mentioned in \cite{Santos-noflips} already
allowed us to obtain polytopes with a non-connected space of triangulation,
but only after a drastic explosion to dimension
234. 

Thanks to symmetries,  we can even be more precise: the graphs of triangulations (and the toric Hilbert schemes)
  of $A_{26}$ and $A_{50}$ have at
  least 17 and 13 connected components
  respectively. In $A_{50}$, each of these 13
  components contains at least $3^{48}$ unimodular triangulations.
The corresponding connected components in the toric Hilbert scheme have dimension
at least 96, in contrast to the fact that the
coherent irreducible component has dimension $50-5-1=44$
(Corollary \ref{coro.quantitative}).

The point set $A_{50}$ has an additional feature: it equals the set of lattice points in 
a lattice polytope $Q$. Point sets of this type play an even more important role in
toric geometry than general lattice point sets. Actually, because 
the construction of \cite{Santos-noflips} did not have this property
the following assertion of V. Alexeev  
\cite[p.705, question 6]{Alexeev} was not completely true: 
``A recent example of F. Santos shows that the moduli 
space $M_Q$  is in general not connected, 
where $Q$ is a lattice polytope". The assertion is true, however,
for the polytope $Q=\conv(A_{50})$:

\begin{theoremintro}
\label{thm:alexeev}
The moduli space of all stable semiabelic toric pairs of type bounded by the product
$Q$ of a segment and a 24-cell (w.r.t. the lattice generated by the vertex set of $Q$)
has at least 13 connected components, each of dimension at least 96 and each with at least
$3^{48}$ torus fixed points.
\end{theoremintro}

\begin{proof}
It follows from Corollary 4.9 of \cite{Santos-noflips} and the discussion preceding it 
that the moduli space in the statement has the same number of connected components as
the graph of triangulations of the set of lattice points in $Q$.

The dimension bound follows the ideas detailed in the proof of Corollary \ref{coro.quantitative} 
for the toric Hilbert scheme. There, we show that in any of the 13 connected
components of the graph of triangulations of $A_{50}$ there are 
triangulations refining a certain subdivision $S_0$ whose only
minimal non-simplicial cells are 48 different octahedra. We then argue on the state polytope
of a certain $A_{50}$-graded ideal whose existence is guaranteed by our Lemma \ref{lemma:unimodular}. 

For the present statement the situation is even simpler:
Theorem 2.13.3 in \cite{Alexeev} says that the 
\emph{generalized secondary polytope} 
(cf. \cite[Section 2.12]{Alexeev} or \cite[Section 4]{Santos-noflips})
of any subdivision of $A_{50}$
is the image of a certain stratum in the moduli space of the statement.
For our subdivision $S_0$, the generalized secondary polytope is 96-dimensional
and has $3^{48}$ vertices: it is a product of  48 triangles.
\end{proof}

It is even possible to obtain a lattice point set that 
\emph{simultaneously} satisfies all the properties stated above.
It is obtained from $A_{50}$ via the 
so-called {\em reoriented Lawrence construction}
introduced in \cite[Section 4.4]{Santos-OMtri}.
%
%

\medskip
Even if the present constructions improve considerably those in
\cite{Santos-noflips}, they still leave open some of the problems
mentioned there:
\begin{itemize}
\item Can the graphs of triangulations be non-connected for point
  sets of \emph{dimensions 3 and 4}? Dimension 3 is specially
  important for 
  Computational Geometry and Engineering applications. 
In dimension 2 the graph of triangulations is easily proved
  to be connected.
\item Can the graphs of triangulations be non-connected for point sets
  \emph{in general position}?  
  General position (i.e., no $d+2$ points lie in an
  affine hyperplane) is interesting in applied areas. Also, 
 a disconnected graph of triangulations for a point set in
general position would imply that the refinement poset of subdivisions 
 of either the whole point set or a proper subset of it is disconnected too.
Connectivity of this poset is still open. Its study is sometimes referred to
as the \emph{generalized Baues problem} for triangulations. 
See \cite{Reiner-survey} or the introduction to \cite{Santos-noflips} for more details.
\item Can the graphs of triangulations be non-connected for \emph{Lawrence
  polytopes}? This would provide realizable oriented matroids with
  non-connected extension space. See more information on this, 
  for example,
  in the introduction to \cite{Santos-noflips}.
\end{itemize}

\section{Triangulations and flips.}
\label{sec.prelim}
\label{sec.subcomplexes}

\subsection{Triangulations and flips of point sets.}
A triangulation of a finite point set $A\subset\reals^d$ is
a geometric simplicial complex with vertex
set contained in $A$ and which covers the convex hull of $A$. 

Geometric bistellar flips are local operations which transform one
triangulation of $A$ into another. Essentially, they correspond to switching
between the two triangulations of a minimal affinely dependent subset of
$A$. More precisely: Every minimal dependent subset $Z\subset A$ 
can be divided
in a unique way in two parts $Z^+$ and $Z^-$ whose convex hulls
intersect. 
We call the ordered pair $(Z^+,Z^-)$ a \emph{circuit} of $A$. (This deviates
slightly from the standard oriented matroid
terminology, in which $Z$ is a circuit and $(Z^+,Z^-)$ an oriented circuit).
The only two triangulations of $Z$ are:
\[
T_{+}(Z):=\{ S\subseteq Z : Z^+\not\subseteq S\},
\qquad
T_{-}(Z):=\{ S\subseteq Z : Z^-\not\subseteq S\}.
\]

\begin{defi}
\label{defi.flips}
Let $T$ be a triangulation of $A$ and let $(Z^+,Z^-) \subseteq A$
be a circuit of $A$.
Suppose that the following conditions are satisfied:
\begin{enumerate}
\renewcommand{\labelenumi}{{\rm(\roman{enumi})}}
\item
  The triangulation $T_{+}(Z)$ is a sub-complex of $T$.
\item
  All the maximal simplices of
$T_{+}(Z)$ have the same link $L$ in $T$. In particular,
$T_{+}(Z)*L$ is a sub-complex of $T$.
Here $A*B:=\{S\cup T: S\in A,\ T\in B\}$ is the \emph{join} of two simplicial complexes.

\end{enumerate}
Then, we can obtain a new triangulation $T'$ of $A$
replacing the sub-complex $T_{+}(Z)*L$ of $T$ by the complex $T_{-}(Z)*L$.
This operation is called a {\em geometric  bistellar
operation} or {\em geometric bistellar flip} (or a {\it flip}, for short)
supported on the circuit $(Z^+,Z^-)$.
\end{defi}

This definition is literally taken from \cite{Santos-noflips}. It
originally comes
from \cite[Chapter 7]{GKZbook}, where it is called a {\em modification}.
See also \cite[p.287]{Lee-survey}.
The following arguments will help convince the reader that it is the right
concept of minimal or elementary change between triangulations:
\begin{enumerate}
\item There is a certain subset of all triangulations of $A$, called {\em
    regular} or {\em coherent}, which are in bijection with 
  the vertex set of a polytope of
  dimension $|A|-d-1$ (the \emph{secondary polytope} of $A$). The edges
of this
  polytope are in bijection with flips between regular
  triangulations \cite{GKZbook}. 

\item Triangulations are the minimal elements in the poset of polyhedral
  subdivisions of $A$, with the partial order given by refinement.
 Flips are the ``next to minimal'' elements, in the following
  well-defined sense: any subdivision whose
only proper refinements are triangulations has exactly two proper refinements
and they are triangulations related by a  flip. Conversely, every flip
arises in this way 
\cite[Corollary 4.5 and Proposition 5.3]{Santos-refine}.
\end{enumerate}

\subsection{Locally acyclic orientations}

Let $A$ be a point set in $\reals^d$. Let $I=[0,1]\subset \reals$.
If $T$ is a triangulation of $A$, we abbreviate as $T\times I$ the 
polyhedral subdivision of $A\times \{0,1\}$ into prisms $\sigma\times I$,
$\sigma\in T$. 
We are interested in studying the triangulations of $A\times \{0,1\}$ that
refine a given such subdivision.

For the refining process we
need to understand triangulations of the prism 
$\Delta^d\times I$, where $\Delta^d$ is a simplex of
dimension $d$. The following description
appears in \cite{Trian-Book} and \cite{GKZbook}; see also \cite[Section
3]{Santos-fewflips}. It can be rephrased as saying that all triangulations of 
$\Delta^d\times I$ are {\em staircase triangulations}.

\begin{proposition}
Let the vertices of $\Delta^d\times I$ be labeled
$\{a_1,\dots,a_{d+1},b_1,\dots$, $b_{d+1}\}$ so that 
the $a_i$'s are the vertices of the facet $\Delta^d\times\{0\}$ and each
$b_i$ is the vertex corresponding to $a_i$ in the opposite facet
$\Delta^d\times \{1\}$. Then:
\begin{enumerate}
\item There is a bijection between triangulations of $\Delta^d\times
I$ and linear orderings (permutations)
of the numbers $\{1,\dots,d+1\}$. The
triangulation corresponding to the ordering $(s_1,\dots ,
s_{d+1})$ has the following $d+1$ maximal simplices:
\[
\left\{
  \{a_{s_{1}}, \dots, a_{s_{i}},b_{s_{i}} \dots, b_{s_{d+1}}\}
\quad:\quad i=1,\dots,d+1
\right\}
.
\]
\item Two triangulations 
of $\Delta^d\times I$ differ by a bistellar flip if and only if the
corresponding orderings differ by a transposition of a pair of
consecutive elements.
\end{enumerate}
\end{proposition}


In particular, we get the following result, where a \emph{locally acyclic
  orientation} of the 1-skeleton of a simplicial complex is an
  orientation of all its edges which is acyclic on every simplex.

\begin{proposition}
\label{prop.locallyacyclic}
The triangulations of $A\times \{0,1\}$ which can be
obtained by refining the product $T\times I$  are in bijection with
the locally acyclic orientations 
of the 1-skeleton of $T$. Flips between such
triangulations correspond to reversal of single edges.
\end{proposition}

\begin{proof}
A locally acyclic orientation induces a linear ordering on every simplex
$\sigma$ of $T$: $i < j$ iff there is an arrow from $i$ to $j$.
We can use this to triangulate $\sigma\times I$. The
triangulations so obtained agree on common faces of any two prisms because
the linear orderings agree on the intersection of any two simplices.
Conversely, a refinement of $T\times I$ triangulates in particular each prism
$\sigma\times I$,
hence it induces a linear ordering on the vertices of every simplex 
$\sigma$ of $T$.
\end{proof}

The triangulation of $A\times\{0,1\}$ obtained from a certain locally acyclic
orientation of $T$ is characterized by the following property: for every
edge $\{v,w\}$ of $T$, directed from $v$ to $w$, the triangulation uses the
diagonal $\{(v,0),(w,1)\}$ in the quadrilateral
$\{(v,0),(w,0),(v,1),(w,1)\}$ of $T\times I$.

We will be specially interested in locally acyclic orientations of $T$
without reversible edges, meaning that every single-edge reversal
creates a cycle in a simplex. The smallest one we
know of has 11 vertices and dimension 3. A slightly bigger example, with
15 vertices, can be obtained as a Schl\"egel diagram of the boundary of the
example considered in Remark 3.4 of \cite{Santos-noflips}.

Unfortunately, Proposition \ref{prop.locallyacyclic}
does not imply that locally acyclic
orientations of $T$ produce refinements of
$T\times I$ without flips. They only produce refinements of $T\times I$ none
of whose flip neighbors refine $T\times I$.

\subsection{Freezing sub-complexes in triangulations}
All we have mentioned so far is essentially present in
\cite{Santos-fewflips}. The new
ingredient in this paper is that we focus on restrictions of
triangulations to sub-complexes of faces. 
Let $F$ be a face of the polytope $\conv(A)$. 
It is obvious that every triangulation of $A$ restricts to a 
triangulation of $F\cap A$. Our next statement says that the same
happens for
flips:

\begin{proposition}
\label{prop.flipfaces}
If $F$ is a face of $\conv(A)$ and $T$ and $T'$ are triangulations of
$A$ differing by a flip, then $T$ and $T'$ restricted to $F$ either
coincide or differ by a flip on a circuit contained in $F$.
\end{proposition}

\begin{proof}
The only simplices of a triangulation that are
removed by a flip are those containing the negative 
part of the circuit $Z$ in which the flip is supported.
For $T$ restricted to $F$ to be affected by the flip it is
necessary that $Z^-\subseteq F$ and, being a face, then $Z^+\subseteq F$
too. Hence, the circuit is contained in $F$ and the flip restricts to
a flip in $F$.
\end{proof}

More generally,
let $\K$ be a simplicial sub-complex of the face complex of 
$\conv(A)$. With the word  
\emph{simplicial} we do not only mean that 
every $F\in \K$ is a simplex, but also that $F\cap A$ is affinely
independent for every $F\in \K$. That is, that $F\cap A$ is the 
vertex set of $F$. By a triangulation of $\K\times I$ 
we mean any geometric simplicial complex $T$ with vertex
set contained in $A\times \{0,1\}$ satisfying that: (1) every simplex of $T$
is contained in one of the products $F\times I$, $F\in \K$; and (2)
$
\cup_{\sigma\in T} \conv(\sigma)=
\cup_{F\in \K} F\times I.
$

The following result has essentially the same proof as 
Proposition \ref{prop.locallyacyclic}, taking into account
Proposition \ref{prop.flipfaces} and the following observation:
for every face $F$ in $\K$, $F\times I$ is a face of  $\conv(A)\times I$.
In particular, every
triangulation of $\K\times I$ induces a triangulation of $F\times I$. 
\begin{theorem}
\label{thm.map}
Let $\K$ be a simplicial subcomplex of the face complex of $\conv(A)$,
for a finite point set $A$.
Triangulations of $\K\times I$ are in bijection with locally acyclic
orientations of the 1-skeleton of $\K$. Flips between triangulations
correspond to locally acyclic orientations differing on the reversal of a 
single edge.
\noproof
\end{theorem}


\begin{corollary}
\label{coro.map}
In the above conditions, 
let $A'$ be any point set in $\reals^{d+1}$ containing $A\times \{0,1\}$
such that for every $F\in \K$ the following two conditions hold:
\begin{itemize}
\item $F\times I$ is still a face of $\conv(A')$, and
\item $F\times I$ contains no point of $A'$ other than its vertex set.
\end{itemize}
Then, ``restriction to $\K\times I$'' induces a simplicial (in particular,
  continuous) map from the graph of triangulations of $A'$ to the graph of
  locally acyclic orientations of $\K$. (Edges in the latter are 
  single-edge reversals).

In particular, if $\K$ has a locally acyclic orientation without reversible
edges and the corresponding triangulation of $\K\times I$ can be
extended to a
triangulation of $A'$, then the graph of triangulations of $A'$ is not
connected.
\end{corollary}

\begin{proof}
The conditions on $A'$ imply that restriction of triangulations of $A'$
to $\K\times I$ is a
well-defined operation and, by Theorem \ref{thm.map}, the restricted
triangulations can be considered elements in the graph of locally acyclic
orientations of $\K$. The map is simplicial by Proposition
\ref{prop.flipfaces}. 

The last sentence in the statement holds because a triangulation
extending (the triangulation of $\K\times I$ corresponding to)
a locally acyclic orientation without reversible edges will not be
connected by flips to a triangulation extending any other locally acyclic
orientation. For example, one extending
any globally acyclic orientation of edges of $\K$, which clearly has
reversible edges and can be extended to a lexicographic 
triangulation of $A'$.
\end{proof}

\subsection{Unimodular triangulations and the toric Hilbert scheme.}
\label{sub.hilbert}

Let $A=\{a_1,\dots,a_n\}\subset\rationals^d$ be a rational point set. 
We transform the \emph{point} configuration
$A$ into a \emph{vector} configuration
$\A=\{\aa_1,\dots,\aa_n\}\subset\rationals^{d+1}$, 
by choosing a positive integer $l_i$ for each $i=1,\dots,n$ and
letting $\aa_i=(l_i a_i,l_i)\in\integers^{d+1}$. We assume, without loss of
generality, that $\A$ has only integer entries.
The standard choice of scaling factors $l_i$ is the
\emph{homogeneous} one, with $l_i=1$ for every $i$, but in Section
\ref{sec.26points} we need a different one.

As detailed in 
\cite[Chapter 10]{Sturmfels-book}, 
we use $\A$ to define a $(d-1)$-dimensional
multi-grading on the polynomial ring
$\field[x_1,\dots,x_n]$, assigning multi-degree $\aa_i$
to the variable $x_i$.
Ideals $I\subset\field[x_1,\dots,x_n]$ that are
homogeneous with respect to this grading have 
a well-defined Hilbert function
\[
\begin{matrix}
\naturals\A & \longrightarrow & \naturals \cr
\bb            &  \longmapsto    & \dim_{\field}I_\bb\cr
\end{matrix}
\]
where $\naturals$ is the set of non-negative integers, $\naturals \A$ is the
semigroup of non-negative integer combinations of $\A$ and for each
$\bb\in\naturals \A$, $I_\bb$ is the part of degree $\bb$ of $I$.

The most natural $\A$-homogeneous ideal is the toric ideal 
$I_\A$ of $\A$, generated by the binomials 
\[
\{x_1^{\lambda_1}\dots x_n^{\lambda_n} - x_1^{\mu_1}\dots x_n^{\mu_n} :
  \lambda,\mu\in \naturals^n, \sum_{i=1}^n (\lambda_i-\mu_i) \aa_i = 0
\}.
\]
For every $\bb\in\naturals \A$, $(I_\A)_\bb$ has
codimension one in $(\field[x_1,\dots,x_n])_\bb$. This characterizes the
Hilbert function of $I_\bb$.

\begin{defi}
\label{defi.agraded}
An $\A$-homogeneous
ideal $I\subset \field[x_1,\dots,x_n]$ is called \emph{$\A$-graded}
if it has the same Hilbert function as the toric ideal.
\end{defi}

Most of the literature on $\A$-graded ideals, starting with
\cite{Sturmfels-book},  assume $\A$ (hence $A$) to have
non-negative entries. In our context this is not important: 
if $A$ has negative entries, then
$\A$ can be mapped to a non-negative configuration by a unimodular
transformation, and  the concept of 
$\A$-graded ideal is invariant under such transformations. 
It is not invariant, however,
under change of choice of the scalars $l_i$. 

$\A$-graded ideals include all the
initial ideals (more generally, all the toric deformations) of $I_\A$. 
The toric Hilbert scheme, as introduced by Peeva and Stillman
\cite{PeeSti}, is
the set of all $\A$-graded ideals with a suitable algebraic
structure defined by some determinantal equations. (An equivalent
 description via
binomial equations appeared in \cite[Section 6]{Sturmfels-agraded}).
See also \cite{MacTho,StStTh}.

Sturmfels \cite[Theorem 10.10]{Sturmfels-book} proved that 
every $\A$-graded ideal $I$ has canonically associated a
polyhedral subdivision $S_I$ of $A$. Observe
that subdivisions of $A$, in the affine setting used in this paper, coincide
with subdivisions of $\A$ in the linear setting 
used in \cite[Chapter 10]{Sturmfels-book}.
If $I$ is monomial, then $S_I$ is a triangulation, whose simplices are the
standard monomials in $\field[ x_1,\dots,x_n]/\Rad(I)$. In other words,
$S_I$ is the triangulation whose Stanley-Reisner ideal
equals the radical of $I$. 

This produces a map 
$\Phi$ from the toric Hilbert scheme of $\A$ 
to the poset of all polyhedral subdivisions of $A$.
(This map factors via the morphism to the Chow variety constructed in
\cite{HaiStu}, followed by the natural map from the Chow variety to the poset
of subdivisions).
Maclagan and Thomas \cite{MacTho} go further and construct a graph of monomial
$\A$-graded ideals ({\em mono-$\A$-GIs} for short) by suitably defining a
concept of flip between mono-$\A$-GIs with the following properties:

\begin{proposition}
\begin{enumerate}
\item The toric Hilbert scheme is connected if and only if the graph of
  mono-$\A$-GIs is connected.
\item The triangulations 
of $A$ corresponding to neighboring mono-$\A$-GIs either
  coincide or differ by a geometric bistellar flip.
\end{enumerate}
\noproof
\end{proposition}

This does not imply that a non-connected graph of
triangulations provides a non-connected toric Hilbert scheme, because
the map $\Phi$ is in general not surjective \cite[Example 10.13]{Sturmfels-book}, 
\cite{PeeSti}.
However, Maclagan and Thomas make the observation, based on \cite[Theorem
10.14]{Sturmfels-book}, that the image of
$\Phi$ contains all the unimodular triangulations of $\A$. A
\emph{unimodular triangulation} of $\A$ is one in which every maximal
simplex is a basis for the lattice spanned by $\A$. 

\begin{corollary}
\label{coro.mactho}
The toric Hilbert scheme has at least as many connected components as
connected components of the graph of triangulations contain unimodular
triangulations.
\noproof
\end{corollary}

To prove the dimension bound  of Corollary \ref{coro.quantitative}
we need the following  generalization of (one direction of) Theorem 10.14 of
\cite{Sturmfels-book}:

\begin{lemma}
\label{lemma:unimodular}
Let $S$ be a subdivision of $\A$ with the property that every cell $\sigma\in S$
can be covered by unimodular simplices with vertex set contained in $\sigma$. Then,
there is an $\A$-graded ideal of the form $\cap_{\sigma\in S} J_\sigma$, where 
each $J_\sigma$ is torus isomorphic to the toric ideal of the set $\sigma$.

In particular, $S$ is in the image of the map $\Phi$.
\end{lemma}

\begin{proof}
The proof follows the ideas in \cite[Theorem 10.14]{Sturmfels-book}.
For each cell $\sigma\in S$, let $I_\sigma$ denote the usual toric ideal of the vertex set of $\sigma$ and 
define
\[
J_\sigma = I_\sigma + \langle x_j : j \not\in \sigma \rangle.
\]
We claim that $I_S:=\cap_{\sigma\in S} J_\sigma$ is indeed $\A$-graded.

 Let $\bb\in \naturals \A$, and let $\sigma$ be a cell of $S$ whose cone 
contains $\bb$. The claim follows from the following three assertions:
(1) all monomials of degree $\bb$ with support not in $\sigma$ are in $I_S$, because
none of them has support contained in a cell of $S$; hence they are in the monomial part of every 
component of $I_S$.
(2) all monomials of degree $\bb$ with support in $\sigma$ are equal modulo $I_S$,
because $I_S$ contains $I_\sigma$.
(3) there are monomials of degree $\bb$ not in $I_S$: by our covering hypothesis, 
there is a unimodular simplex $\tau$ contained in
$\sigma$ and whose cone contains $\bb$. Consider the unique positive combination 
$\sum_{i\in\tau} \lambda_i \aa_i$
that gives $\bb$. This combination is integer since $\tau$ is unimodular. The corresponding
monomial $\Pi_{i\in \sigma} x_i^{\lambda_i}$ is clearly not in  
$I_\sigma + \langle x_j : j \not\in \sigma \rangle=J_\sigma$, hence not in
$I_S$.

That $\Phi(I_S)=S$ is essentially the definition of $\Phi$; see \cite[Theorem 10.10]{Sturmfels-book}.
\end{proof}

\section{A construction with 50 points, based on the 24-cell}
\label{sec.50points}
The 24-cell is one of the six regular 4-dimensional polytopes.
Its 24 vertices are:
\begin{itemize}
\item The 8 points $\pm 2 e_i$.
\item The 16 points $(\pm 1,\pm 1,\pm 1,\pm 1)$.
\end{itemize}
Let $A$ consist of these 24 vertices. Let $O$ be the origin and
let $\K$ be the 2-skeleton of the 24-cell, 
which consists of 96 triangles and 96 edges. Let 
\[
A_{50}:=(A\cup\{O\})\times\{0,1\}.
\]

\begin{theorem}
\label{thm.24cell}
The 5-dimensional point configuration with 50 elements 
$A_{50}$ has a non-connected space
of triangulations and its homogenized version $\A_{50}=A_{50}\times \{1\}$ has
a non-connected toric Hilbert scheme.
\end{theorem}

\begin{proof}
We show below that the 2-skeleton $\K$ of the 24-cell can be
given a locally acyclic orientation with no reversible edges
(Lemma \ref{lemma.orient24}) and that the triangulation
of $\K\times I$ corresponding to this orientation 
can be extended to a unimodular triangulation $T'$ of $A_{50}$
(Lemma \ref{lemma.extend24}).
Since $A_{50}$ satisfies the conditions required for $A'$ in
 Corollary \ref{coro.map}, its graph of triangulations
 is not connected. Unimodularity of $T'$
implies that the toric Hilbert
scheme of $\A_{50}:=A_{50}\times \{1\}\subset\reals^6$ is not
connected either, via Corollary \ref{coro.mactho}.
\end{proof}

In the rest of this section we fill in the details in this proof and prove more precise
quantitative results, and that the graph of triangulations of $A_{50}$
remains not connected when its only two non-vertices $(O,0)$ and $(O,1)$ are
moved into convex position.

\medskip

The facets of the 24-cell are 24 octahedra. One of them
is the octahedron with vertices $(2,0,0,0)$,
$(0,2,0,0)$, $(1,1,1,1)$, $(1,1,-1,-1)$, $(1,1,1,-1)$ and $(1,1,-1,1)$.
We denote it $F_{1,1,0,0}$.
From this we get three other octahedra $F_{-1,1,0,0}$, $F_{-1,-1,0,0}$
and $F_{1,-1,0,0}$ by the rotation of order 4 on the first two coordinates.
And from these four we get the rest of the octahedra by permuting
coordinates. 

The subindices in each octahedron give the coordinates of
its centroid. The 24-cell is self-polar: These 24
centroids 
are the vertices of another regular (and smaller) 24-cell. 
Our  coordinates are chosen to  highlight the
symmetries of the 24-cell that we are interested in.

\begin{figure}[ht]
  \epsfysize = 5 cm
  \leavevmode
  \epsfbox{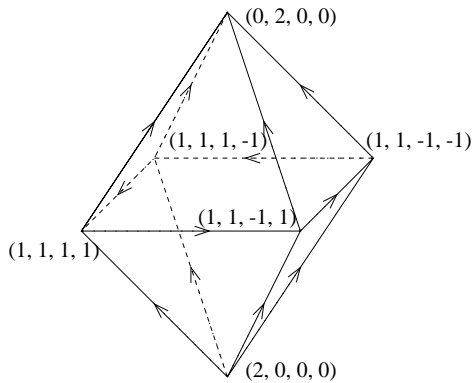}
\caption{A locally acyclic orientation of the 1-skeleton of the octahedron
  $F_{1,1,0,0}$,
  with only 4 reversible edges.}
\label{fig.octahedron}
\end{figure}

We orient the edges in $F_{1,1,0,0}$ with a source at $(2,0,0,0)$, a
sink at $(0,2,0,0)$ and a cycle of length 4 on the equatorial square
$(1,1,1,1)\to(1,1,-1,1)\to(1,1,-1,-1)\to(1,1,1,-1)$. 
See Figure \ref{fig.octahedron}.
Observe that among the 12
edges of this octahedron only the equatorial four 
can be reversed without creating a local cycle.
We let the other 84 edges of the 24-cell
be oriented by the action of the affine group $G$ of order 32
generated by the exchange of first two and last two
coordinates and the rotation of order 4 on the plane of the first two
coordinates (or of the last two coordinates).

\begin{lemma}
\label{lemma.orient24}
This orientation of the 1-skeleton of $\K$ is well-defined, locally
acyclic, and has no reversible edges.
\end{lemma}

\begin{proof}
The orientation is well-defined because edges of $F_{1,1,0,0}$
in the same orbit under $G$ receive orientations compatible with the action.
Indeed, the only symmetries in $G$ sending some edge of
$F_{1,1,0,0}$ to another are the rotations on the last two coordinates
(rotation around the vertical axis, in Figure \ref{fig.octahedron})
and they are compatible with the orientation given.

Under $G$, there are two orbits of octahedra in the 24-cell, 3 orbits
of triangles and edges, and two orbits of vertices. The
following are representatives of the three orbits of triangles, with
their orientations indicated:
\begin{itemize}
\item $(2,0,0,0) \to(1,1,1,-1)\to (1,1,1,1)$,
\item $(1,1,1,-1)\to (1,1,1,1)\to (0,2,0,0)$, and
\item $(1,1,1,-1)\to (0,0,2,0)\to (1,1,1,1)$.
\end{itemize}
This proves that the orientation is locally acyclic. The three orbits
of edges have as representatives $(2,0,0,0) \to(1,1,1,1)$,
$(1,1,1,-1)\to (0,2,0,0)$, and $(1,1,1,-1)\to (1,1,1,1)$. That is, the
non-reversible edges of the three stated representatives of triangles. This
proves that the orientation has no reversible edges.
\end{proof}

Actually, the symmetry group of the locally acyclic
orientation we are using is
larger than $G$. We leave it to the interested reader to check that it
is an orientation-preserving group of order 96
and acts transitively over the 96 edges and over
the 96 triangles in the 24-cell. It is generated by $G$ and, for
example, the
following orthogonal transformation:
\[
\left(
\begin{matrix}
x_1 \cr x_2 \cr x_3 \cr x_4
\end{matrix}
\right)
\longmapsto
\frac{1}{2}
\left(
\begin{matrix}
x_1 - x_2 + x_3 + x_4\cr
x_1 + x_2 - x_3 + x_4\cr
x_1 + x_2 + x_3 - x_4\cr
-x_1 + x_2 + x_3 + x_4\cr
\end{matrix}
\right)
\]

\medskip


\begin{lemma}
\label{lemma.extend24}
The triangulation  $T$ of $\K\times I$ corresponding 
via Theorem \ref{thm.map} to this locally acyclic orientation 
of $\K$
can be extended to a unimodular triangulation $T'$ of $A_{50}$.
\end{lemma}

In the above statement, a triangulation of a lattice point set 
$B\subset\reals^d$ is called
unimodular if every maximal simplex is an affine
basis for the affine lattice spanned by $B$.
Equivalently, if it is unimodular as a triangulation of the vector set
$\mathcal B$ obtained from $B$ with a homogeneous choice of scaling lengths.

\begin{proof}
Let $T_{1,1,0,0}$ be the triangulation of the octahedron
$F_{1,1,0,0}$ that uses the 
axis $\{(2,0,0,0),(0,2,0,0)\}$, oriented from 
$(2,0,0,0)$ to $(0,2,0,0)$.
This triangulation extends in a locally acyclic way the
orientation of edges of $F_{1,1,0,0}$.

Let $T_0$ be the triangulation of $A\cup\{O\}$ obtained by first
replicating $T_{1,1,0,0}$ to the other octahedra by the action of $G$ and
then coning the triangulated boundary of the 24-cell to the centroid $O$.
We also extend the orientation to all the boundary via $G$, 
and orient all the
edges incident to $O$ with a source at $O$. This gives a locally acyclic
orientation of the 1-skeleton of $T_A$ and, hence, a triangulation $T'$ of
$A_{50}$ which refines $T_0\times I$. Since the orientation of $T_0$
extends the one we had in $\K$, $T'$ extends what we had in $\K\times I$.

Only unimodularity remains to be
checked. Since a prism over a unimodular simplex is totally
unimodular (meaning that all its triangulations are unimodular), every
refinement of $T_0\times I$ will be unimodular as long as $T_0$ itself is
unimodular. By symmetry,
we only need to check that one of the 96 maximal simplices of $T_0$
is unimodular. 
 For example, take the 4-simplex with
vertex set $(0,0,0,0)$, $(2,0,0,0)$, $(0,2,0,0)$, $(1,1,1,1)$,
$(1,1,1,-1)$. This is unimodular (in the sub-lattice spanned by
$A\cup\{0\}$) 
because
\[
\left\vert
\begin{matrix}
  0 &  2 &  0 &  1 &  1 \cr
  0 &  0 &  2 &  1 &  1 \cr
  0 &  0 &  0 &  1 &  1 \cr
  0 &  0 &  0 &  1 & -1 \cr
  1 &  1 &  1 &  1 &  1 \cr
\end{matrix}
\right\vert
=- 8
\]
and the sub-lattice spanned by
$A\cup\{O\}$, given by the conditions ``all coordinates have the same
parity'', has index 8
in the integer lattice.
\end{proof}


That finishes the proof of Theorem~\ref{thm.24cell}. 
We now prove the quantitative results mentioned after
Theorem \ref{thm.main} and that the graph of triangulations of $A_{50}$
remains unimodular in a convex-position version of the point set.

\begin{corollary}
The graph of triangulations of $A_{50}$ and the toric Hilbert scheme of its
homogenized version $\A_{50}=A_{50}\times \{1\}$ have at least 13 connected
components, each with at least $3^{48}$ unimodular 
triangulations/monomial ideals. In the toric Hilbert scheme, each of these 13
connected components
has dimension at least 96.
\label{coro.quantitative}
\end{corollary}

\begin{proof}
Observe in Figure \ref{fig.octahedron} that, out of the 48 symmetries
of the octahedron $F_{1,1,0,0}$, only the four generated by the
rotation around the vertical axis leave our locally acyclic
orientation invariant. Hence, there are 12 ways of constructing a
locally acyclic orientation with exactly the same properties. (This
agrees with our claim that our locally acyclic orientation has 96
symmetries, since the symmetry group of the 24-cell has $24\times 48
=12\times 96$ elements). Each of these 12 orientations is a different
isolated element in the graph of locally acyclic orientations of $\K$,
hence it produces a different connected component in the graph of
triangulations of $A_{50}$. To these 12 we have to add the regular
component of the graph, obtained for example by starting with a
globally acyclic orientation of $\K$.

To prove the number $3^{48}$, observe the following: 
let $F$ be one of the 24 octahedra in the 24-cell.
Let $v$ and $w$ be the source and sink of its orientation.
In the triangulation $T_0$ of the proof of Lemma
\ref{lemma.extend24}, $F$ is triangulated into
four 3-simplices which are joined to the origin $O$.
In particular, the two flips of the triangulation of $F$ are
still flips in $T_0$. 

Moreover, these four 4-simplices of $T_0$ have the
same source $O$ and sink $w$. Hence,
the copies of them in $(A\cup O)\times \{0\}$ are joined
in $T'$ to the same vertex of $(A\cup O)\times \{1\}$
(namely $(w,1)$) and the copies of them 
in $(A\cup O)\times \{1\}$ are joined
to the same vertex of $(A\cup O)\times \{0\}$
(namely $(O,0)$). 
This implies that the two flips in $F$ produce four
flips in $T'$: Two with supporting circuit contained in 
the octahedron $F\times\{0\}$ and another two with supporting 
circuit in
$F\times\{1\}$. Moreover, flips in different octahedra, out of
the total of 48, are independent. Hence, from $T'$ we can at least reach the
$3^{48}$ (including $T'$) triangulations obtained by these flips, which are all
unimodular. 

For the dimension bound, let $S_0$ be the subdivision of $A_{50}$
obtained from $T_0$ by making each of the 48 octahedra (and the two
vertices coned to each of them) a single cell. Any refinement of this
subdivision is unimodular, because the three triangulations of each
octahedron are unimodular. By Lemma \ref{lemma:unimodular}, there is
an $\A_{50}$-graded ideal $I_0$ associated to that subdivision. Since
$I_0$ is the sum of 48 independent copies of the toric ideal of an
octahedron, the state polytope of $I_0$ (cf. \cite{Sturmfels-book}) is
the product of 48 copies of a triangle, hence it has dimension 96. And
the toric variety of the state polytope is immersed (modulo normalization)
in the toric
Hilbert scheme because every toric deformation of an $\A$-graded ideal
is $\A$-graded as well.
\end{proof}

We can even be more precise; the subgraph induced by the $3^{48}$ triangulations/ideals
mentioned in the statement of Corollary \ref{coro.quantitative} is the 1-skeleton of the
polytope $(\Delta^2)^{48}$, where $\Delta^2$ is a triangle.

\begin{proposition}
\label{prop.50convex}
Let $A'_{50}$ be the point set obtained from $A_{50}$ by moving the
points 
$(O,0)$ and $(O,1)$ to $(O,\alpha)$ and $(O,\beta)$, for any $\alpha<0$ and
$\beta>1$. Then, $A'_{50}$ is in convex position and its graph of
triangulations still has at least 13 components, each with at least $3^{48}$ triangulations.
\end{proposition}

\begin{proof}
$A'_{50}$ is still a point set containing $A\times \{0,1\}$,
and $\K$ satisfies the two hypotheses in Corollary \ref{coro.map}. The only
thing to prove is that our triangulation $T$ of $\K\times I$ extends to a
triangulation of $A'_{50}$. Actually, the following is true: the same
triangulation $T'$ of Lemma \ref{lemma.extend24}, considered in $A'_{50}$
with the substitution $(O,0)\mapsto (O,\alpha)$
and $(O,1)\mapsto (O,\beta)$ is still a triangulation of $A'_{50}$.
This holds because the facets $A\times \{0\}$ and
$A\times \{1\}$ of $A_{50}$
are centrally triangulated in $T'$ and the
perturbed points 
$(O,\alpha)$ and $(O,\beta)$ lie beyond these facets of $A_{50}$
(with the standard meaning of beyond
in polytope theory: a point lies beyond a facet $F$
of the polytope $P$ if only that facet is visible from the point).
\end{proof}

\begin{remark}
\label{rem.convexscheme}
We cannot say anything about the toric Hilbert scheme of $A'_{50}$,
even if we assume $\alpha$ and $\beta$ to be rational or even integer,
because our triangulations are no longer unimodular in $A'_{50}$.  
\end{remark}

%

\section{A construction with 26 points, based on the cross-polytope}
\label{sec.26points}

In this section we let $A\subset \reals^4$ denote the vertex set of a
regular cross-polytope, with centroid $O$:
\[
A:=\{\pm e_1, \pm e_2, \pm e_3, \pm e_4\}\subset\reals^4,
\qquad\qquad O=(0,0,0,0)\in\reals^4.
\]
The faces of the cross-polytope are 16 tetrahedra, 32 triangles, 24
edges and 8 vertices. Our complex $\K$ will contain all the 24 edges,
but only 24 of the
32 triangles. The point set
$A_{26}\subset\reals^4$ will contain the 18 points in 
\[
\left(A\cup\{O\}\right)\times \{0,1\}
\]
plus 8 points of the form $(p,1/2)$ where $p$ (essentially) runs over
the centroids of the eight triangles missing in $\K$. 
See the definitions of $\K$ and $A$ below.

We will show that $\K$ can be
given a locally acyclic orientation with no reversible edges, and that
this orientation can be extended to a triangulation of $A_{26}$. Also, that
this triangulation is unimodular when considered in a 
vector configuration $\A_{26}$ obtained from $A_{26}$ as in Section
\ref{sub.hilbert}, but with a non-homogeneous choice of scaling factors.
These facts imply that neither the graph of triangulations 
of $A_{26}$ nor the toric Hilbert
scheme of $\A_{26}$ are connected. 

Our construction is symmetric under the group $G$ of order six
generated by the central symmetry and the rotation 
which cyclically permutes the first three coordinates:
\[
(x_1,x_2,x_3,x_4)\mapsto (x_2,x_3,x_1,x_4).
\]
This group produces 4 orbits of edges, with
six members each, and 6 orbits of triangles, one with two elements and five with six elements,
We show them as columns in  Tables \ref{table.edges24} and \ref{table.triangles24}. 
To save space we write $\e_i$ meaning $-e_i$.

\begin{table}[ht]
\begin{tabular}{cccc}
$\{e_1,e_2\}$ & $\{\e_2,e_1\}$ & $\{e_4,e_1\}$ & $\{\e_1,e_4\}$ \cr
$\{e_2,e_3\}$ & $\{\e_3,e_2\}$ & $\{e_4,e_2\}$ & $\{\e_2,e_4\}$ \cr
$\{e_3,e_1\}$ & $\{\e_1,e_3\}$ & $\{e_4,e_3\}$ & $\{\e_3,e_4\}$ \cr
$\{\e_1,\e_2\}$ & $\{e_2,\e_1\}$ & $\{\e_4,\e_1\}$ & $\{e_1,\e_4\}$ \cr
$\{\e_2,\e_3\}$ & $\{e_3,\e_2\}$ & $\{\e_4,\e_2\}$ & $\{e_2,\e_4\}$ \cr
$\{\e_3,\e_1\}$\ & $\{e_1,\e_3\}$\ & $\{\e_4,\e_3\}$\ & $\{e_3,\e_4\}$\ 
\cr
\end{tabular}
\medskip
\caption{The 24 edges of the cross-polytope, divided into four orbits.}
\label{table.edges24}
\end{table}

\begin{table}[ht]
\begin{tabular}{cccccc}
                                  & $\{ e_1,\e_4,\e_2\}$ & $\{ e_1,\e_3, e_2\}$ & $\{\e_3, e_4, e_2\}$ & $\{ e_4, e_1, e_2\}$ & $\{ e_1, e_2,\e_4\}$ \cr
                                  & $\{ e_2,\e_4,\e_3\}$ & $\{ e_2,\e_1, e_3\}$ & $\{\e_1, e_4, e_3\}$ & $\{ e_4, e_2, e_3\}$ & $\{ e_2, e_3,\e_4\}$ \cr
$\{ e_1, e_2, e_3\}$ & $\{ e_3,\e_4,\e_1\}$ & $\{ e_3,\e_2, e_1\}$ & $\{\e_2, e_4, e_1\}$ & $\{ e_4, e_3, e_1\}$ & $\{ e_3, e_1,\e_4\}$ \cr
$\{\e_1,\e_2,\e_3\}$ & $\{\e_1, e_4, e_2\}$ & $\{\e_1, e_3,\e_2\}$ & $\{ e_3,\e_4,\e_2\}$ & $\{\e_4,\e_1,\e_2\}$ & $\{\e_1,\e_2, e_4\}$ \cr
                                  & $\{\e_2, e_4, e_3\}$ & $\{\e_2, e_1,\e_3\}$ & $\{ e_1,\e_4,\e_3\}$ & $\{\e_4,\e_2,\e_3\}$ & $\{\e_2,\e_3, e_4\}$ \cr
                                  & $\{\e_3, e_4, e_1\}$ & $\{\e_3, e_2,\e_1\}$ & $\{ e_2,\e_4,\e_1\}$ & $\{\e_4,\e_3,\e_1\}$ & $\{\e_3,\e_1, e_4\}$ \cr

\end{tabular}
\medskip
\caption{The 32 triangles of the cross-polytope, divided into six orbits.}
\label{table.triangles24}
\end{table}

%

\begin{figure}[ht]
  \epsfysize = 5 cm
  \leavevmode
  \epsfbox{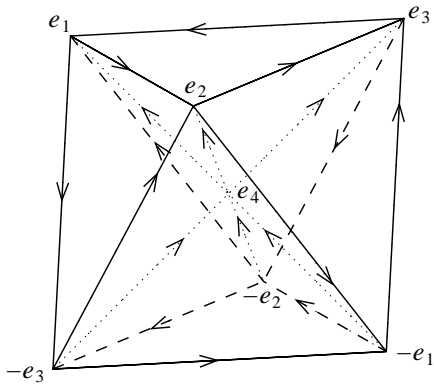}
\caption{A locally acyclic orientation, without reversible edges,
of a sub-complex of the
  2-skeleton of a cross-polytope.}
\label{fig.crosspoly}
\end{figure}

We orient edges in the way 
implicit in Table \ref{table.edges24}. That is, each edge is
oriented from the first vertex listed to the second
(this is clearly compatible with the action of $G$).
Figure \ref{fig.crosspoly} may help understand the construction.
It shows the star of $e_4$ in the boundary of the cross-polytope,
which is an octahedron centrally triangulated into eight tetrahedra.
This is half of the boundary of the cross-polytope; the other half is obtained by
central symmetry and is not drawn. The edges in the figure have been oriented,
and the reader can verify that only five of the triangles in the figure 
have cyclic orientations; namely, those belonging to the first two orbits
of triangles in Table \ref{table.triangles24}. Triangles in the last four orbits
receive the acyclic orientation corresponding to the ordering
in which their vertices are listed in the table.

Let $\K$ be the union of the last four orbits of triangles in Table \ref{table.triangles24}.

\begin{proposition}
\label{prop.lao24}
The orientation given to the 1-skeleton  is locally acyclic (in $\K$)
and has no reversible edges.
\end{proposition}

\begin{proof}
Check that indeed the last four orbits of triangles in Table
\ref{table.triangles24} receive an acyclic orientation
and that the
non-reversible edges of these four orbits lie respectively in the four
orbits of edges 
(only one representative triangle of each orbit needs to be checked).
\end{proof}

As in Section \ref{sec.50points}, the symmetry group of this locally acyclic
orientation is actually larger than $G$. Let us denote it $\tilde G$.
Apart from $G$, it contains for example the 
simultaneous rotation of order 4 in the planes of the first two and
last two coordinates:
\[
\rho(x_1,x_2,x_3,x_4) = (-x_2,x_1,-x_4,x_3).
\]

\begin{proposition}
\label{prop.symmetry26}
$\tilde G$ is an orientation-preserving 
group of order 24 generated by $G$ and  the above transformation
$\rho$. It acts transitively over the 24 edges and 24 triangles of $\K$,
and over the eight triangles not in $\K$.
\end{proposition}

\begin{proof}
Since $\rho$ sends the edges $\{e_1,e_2\}$, $\{e_2,e_3\}$ and 
$\{e_3,e_1\}$ of the 
first $G$-orbit in Table \ref{table.edges24} to edges in the other three
$G$-orbits, respectively,
 the group generated by $\rho$ and $G$ 
(and hence $\tilde G$) acts transitively over the 24 edges of the
24-cell. Moreover, $\tilde G$ acts with trivial stabilizer on every edge,
because every edge is the
non-reversible edge of a unique triangle of $\K$, and the reflection on the
hyperplane containing that triangle and the origin is not in $\tilde G$
(to check this consider any single triangle of $\K$, for
example 
$\{-e_3,e_4,e_2\}$). 

This proves that $\tilde G$ is generated by $G$ and
$\rho$, that it 
has order 24, and that it acts transitively over edges. 
It also acts transitively on triangles of $\K$ because
every edge is the non-reversible edge of a different triangle.
It acts transitively on the remaining eight triangles because $\rho$ sends the
member $\{-e_1,-e_2,-e_3\}$ of the first $G$-orbit to the member
$\{e_1,-e_2,-e_4\}$ of the second $G$-orbit.
\end{proof}

Let $T$ denote the triangulation of $\K\times I$ corresponding to the
locally acyclic orientation of the 1-skeleton of $\K$ we have described.
It is clear that $T$ cannot be extended to a triangulation of 
$(A\cup\{O\})\times\{0,1\}$, because the prism $\{e_1,e_2,e_3\}\times\{0,1\}$,
for example, has its three quadrilateral faces triangulated in a 
non-extendable way. The same holds for the other 7 triangles not in
$\K$. This is why we need to define the
following point set $A_{26}$ with 26 elements, the first 18 of which are 
$(A\cup\{O\})\times\{0,1\}$. Let $A_{26}$ consist of:

\begin{itemize}
\item The 8 points $\pm a_i:=(\pm e_i, 0)$, $i=1,2,3,4$.
\item The 8 points $\pm b_i:=(\pm e_i, 1)$, $i=1,2,3,4$.
\item The two points $a_0:=(0,0,0,0,0)$ and  $b_0:=(0,0,0,0,1)$.
\item The following eight points:
\begin{center}
\begin{tabular}{lr}
$c_{+,+,+,0}=( \frac{1}{2}, \frac{1}{2}, \frac{1}{2},0,\frac{1}{2})$,\ &
$c_{-,-,-,0}=(-\frac{1}{2},-\frac{1}{2},-\frac{1}{2},0,\frac{1}{2})$, 
\medskip\\
$c_{+,-,0,-}=( \frac{1}{2},-\frac{1}{2},0,-\frac{1}{2},\frac{1}{2})$,\ &
$c_{-,+,0,+}=(-\frac{1}{2}, \frac{1}{2},0, \frac{1}{2},\frac{1}{2})$, 
\medskip\\
$c_{-,0,+,-}=( -\frac{1}{2},0,\frac{1}{2},-\frac{1}{2},\frac{1}{2})$,\ &
$c_{+,0,-,+}=(\frac{1}{2},0, -\frac{1}{2}, \frac{1}{2},\frac{1}{2})$, 
\medskip\\

$c_{0,+,-,-}=(0, \frac{1}{2},-\frac{1}{2},-\frac{1}{2},\frac{1}{2})$, &
$c_{0,-,+,+}=(0,-\frac{1}{2}, \frac{1}{2}, \frac{1}{2},\frac{1}{2})$.
 \medskip\\
\end{tabular}
\end{center}
\end{itemize}
Observe that $A_{26}$ is compatible with the action of $\tilde G$
(times the trivial group acting on the fifth coordinate).

Recall that a point $p$ is said to lie 
\emph{beyond} a certain facet $F$ of a
polytope $P$ if
it lies outside $P$ and $F$ is the only facet visible from $p$. We
generalize this and say that $p$ lies beyond a face $F$ if the facets
visible from $p$ are precisely those containing $F$. We do not allow $p$
to lie in any facet-defining hyperplane.

\begin{lemma}
\label{lemma.beyond}
\begin{enumerate}
\item Each point $c_{*,*,*,*}$ lies beyond a face $\sigma\times I$
of $\conv(A)\times I$, where $\sigma$ is one of the 
8 triangles of $\conv(A)$ missing in $\K$.
\item For each triangle $\sigma\in \K$, $\sigma\times I$
  is still a face in $\conv(A_{26})$ and contains no point of $A_{26}$ other
  than its vertex set. 
\end{enumerate}
\end{lemma}

\begin{proof}
By symmetry, we need to prove part (2) only for one triangle in $\K$ and
part (1) only for  one not in $\K$.

For part (1), let $\sigma=\conv(\{e_1,e_2,e_3\})$. The statement is
equivalent to saying
that $(1/2,1/2,1/2,0)$ lies beyond $\sigma$ in $\conv(A)$. This is easy
to check: the facet-defining half-spaces of $\conv(A)$ are those of the
form $\pm x_1 \pm x_2 \pm x_3 \pm x_4\le 1$, and the only ones not
containing $(1/2,1/2,1/2,0)$ are $+ x_1 + x_2 + x_3 + x_4\le 1$ and $+
x_1 + x_2 + x_3 - x_4\le 1$.

For part (2), let $\sigma=\{e_4,e_1,e_2\}$. We leave it to the reader 
to check that the functional 
\[
x_1 + x_2 -x_3/2 + x_4
\]
on $A_{26}$
is maximized exactly on the six vertices of $\sigma\times \{0,1\}$.
\end{proof}

Let $\A_{26}\subset\reals^6$ be the vector configuration 
obtained from $A_{26}$ with scaling factor
equal to 1 for the points $a_*$ and $b_*$ and equal to 2 for the points
$c_{*,*,*,*}$. In other words, let $\A$ consist of:
\begin{itemize}
\item The 8 vectors $\pm \aa_i:=(\pm e_i, 0, 1)$, $i=1,2,3,4$.
\item The 8 vectors $\pm \bb_i:=(\pm e_i, 1, 1)$, $i=1,2,3,4$.
\item The two vectors $\aa_0:=(0,0,0,0,0,1)$ and 
$\bb_0=(0,0,0,0,1,1)$.
\item The following eight vectors:
\begin{center}
\begin{tabular}{lr}
$\cc_{+,+,+,0}=( 1, 1, 1,0,1,2)$,\ &
$\cc_{-,-,-,0}=(-1,-1,-1,0,1,2)$, 
\medskip\\
$\cc_{+,-,0,-}=( 1,-1,0,-1,1,2)$,\ &
$\cc_{-,+,0,+}=(-1, 1,0, 1,1,2)$, 
\medskip\\
$\cc_{-,0,+,-}=( -1,0,1,-1,1,2)$,\ &
$\cc_{+,0,-,+}=(1,0, -1, 1,1,2)$, 
\medskip\\
$\cc_{0,+,-,-}=(0, 1,-1,-1,1,2)$, &
$\cc_{0,-,+,+}=(0,-1, 1, 1,1,2)$.
 \medskip\\
\end{tabular}
\end{center}
\end{itemize}

\begin{theorem}
\label{thm.26points}
\begin{enumerate}
\item The triangulation $T$ of $\K\times I$ can be extended to a
triangulation of
  $A_{26}$. 
\item This extended triangulation is unimodular when considered in 
 $\A_{26}$.
\end{enumerate}
\end{theorem}

Before proving this, let us show its implications:

\begin{corollary}
\label{coro.26points}
The graph of triangulations of $A_{26}$ and the toric Hilbert scheme of
$\A_{26}$ have at least 17 connected components.
\end{corollary}

\begin{proof}
Condition (2) in Lemma \ref{lemma.beyond} says that 
we can apply Corollary \ref{coro.map} to $\K$ and $A_{26}$. 
By that corollary, the triangulation extending $T$ cannot be
connected by flips to any triangulation with a different restriction to
$\K\times I$, which implies that the graph of triangulations of $A_{26}$ (and, by
unimodularity, 
the toric Hilbert scheme of $\A_{26}$) is not connected.

The number 17 comes from the fact that the symmetry group of our locally
acyclic orientation has order 24 (Proposition \ref{prop.symmetry26}) versus
the $16\times 24$ symmetries of the cross-polytope: There are 15 other
equivalent ways of doing our construction, and any regular triangulation of
$A_{26}$ provides a 17th connected component.
\end{proof}

\begin{proof}[Proof of Theorem \ref{thm.26points}]
Let us consider the extended complex $\K'=\K*O$. The orientation in $\K$ can
be extended trivially to a locally acyclic orientation of $\K'$ by
letting $O$
be a global source. We extend $T$ to $\K'\times
I$ using this orientation. This is compatible with the central
triangulations of the 
two facets $(A\cup\{O\})\times \{0\}$ and $(A\cup\{O\})\times \{1\}$.
Hence, we
have a triangulation, that we denote  $T'$, of
\[
(\K'\times I) \cup (\conv\{A\}\times \{0,1\}).
\]

We claim that the complement of this polyhedral complex in
$\conv(A_{26})$ consists of
eight convex regions, each homeomorphic to a 5-dimensional 
closed half-space and with
one of the eight points $c_{*,*,*,*}$ on its boundary. This claim
implies that we can
extend $T'$ to $A_{26}$ by coning each point $c_{*,*,*,*}$ to the
part of $T'$ visible from it.

To prove the claim, observe that the complement of $\K$ in the boundary
of $\conv(A)$
consists of eight regions (homeomorphic to 3-balls), each being the
union of
the two tetrahedra incident to one of the points $c_{*,*,*,*}$. Hence, the
complement of $\K'$ in $\conv(A)$ consists of eight regions
$R_{*,*,*,*}$, one
for each point $c_{*,*,*,*}$ defined in the following way for $c_{+,+,+,0}$
(with the obvious generalization to the other seven regions). 
$R_{+,+,+,0}$ equals the closed region
\[
\conv(O,e_1,e_2,e_3,e_4)\cup\conv(O,e_1,e_2,e_3,-e_4).
\]
minus the part of its boundary incident to $O$. In particular, 
$R_{+,+,+,0}$ is convex and homeomorphic to a 4-dimensional half-space.

By Lemma \ref{lemma.beyond}  the boundary of $\conv(A_{26})$ can be considered
a stellar subdivision of the boundary of $\conv(A\times I)$, obtained by
pulling the 
centroids of certain faces out to the points $c_{*,*,*,*}$.
 In other words,  the complement of
$(\K\times I) \cup (\conv\{A\}\times \{0,1\})$ in 
$\conv(A_{26})$ consists of eight regions, each being the above region
$R_{*,*,*,*}$ times the open segment $(0,1)$
except it has been slightly pulled out
from its interior point $c_{*,*,*,*}$. This finishes the proof of the claim.

Summing up, we have proved that $T$ extends to the triangulation of $A_{26}$
obtained by letting $\tilde G$ (times the trivial group in the fifth
coordinate) act over the twenty-eight 5-simplices listed in Table 
\ref{table.trian26}. All the simplices are meant to contain
$c_{+,+,+,0}$, which we omit.

\begin{table}
\[
\begin{matrix}
\{a_0,a_4,a_1,a_2, b_2\},\ &\ 
\{a_0,a_4,a_2,a_3, b_3\},\ &\ 
\{a_0,a_4,a_3,a_1, b_1\},\ &\ 
\\ 
\{a_0,a_4,a_1,b_1, b_2\},\ &\ 
\{a_0,a_4,a_2,b_2, b_3\},\ &\ 
\{a_0,a_4,a_3,b_3, b_1\},
\\ 
\{a_0,a_4,b_4,b_1, b_2\},\ &\ 
\{a_0,a_4,b_4,b_2, b_3\},\ &\ 
\{a_0,a_4,b_4,b_3, b_1\},
\\ 
\{a_0,b_0,b_4,b_1, b_2\},\ &\ 
\{a_0,b_0,b_4,b_2, b_3\},\ &\ 
\{a_0,b_0,b_4,b_3, b_1\},
\medskip
\\
\{a_0,a_1,a_2,-a_4,-b_4\},\ &\ 
\{a_0,a_2,a_3,-a_4,-b_4\},\ &\ 
\{a_0,a_3,a_1,-a_4,-b_4\},
\\ 
\{a_0,a_1,a_2,b_2,-b_4\},\ &\ 
\{a_0,a_2,a_3,b_3,-b_4\},\ &\ 
\{a_0,a_3,a_1,b_1,-b_4\},
\\ 
\{a_0,a_1,b_1,b_2,-b_4\},\ &\ 
\{a_0,a_2,b_2,b_3,-b_4\},\ &\ 
\{a_0,a_3,b_3,b_1,-b_4\},
\\ 
\{a_0,b_0,b_1,b_2,-b_4\},\ &\ 
\{a_0,b_0,b_2,b_3,-b_4\},\ &\ 
\{a_0,b_0,b_3,b_1,-b_4\},
\medskip
\\
& \{a_0,a_1,a_2,a_3,a_4\} \\
& \{a_0,a_1,a_2,a_3,-a_4\} \\
& \{b_0,b_1,b_2,b_3,b_4\} \\
& \{b_0,b_1,b_2,b_3,-b_4\} \\
\end{matrix} 
\]
\label{table.trian26}
\caption{The triangulation of $A_{26}$.}
\end{table}

The first six groups of four simplices in the list come from the six
tetrahedra 
\[
\begin{matrix}
\{O,e_4,e_1,e_2\}, &\{O,e_4,e_2,e_3\},& & \{O,e_4,e_3,e_1\},\cr
\{O,e_1,e_2, -e_4\}, &\{O,e_2,e_3, -e_4\}&\hbox{ and} &\{O,e_3,e_1, -e_4\}, 
\end{matrix}
\]
which appear in $\K'$ with their vertices ordered in this way.
The last four simplices come from the central
triangulations of $\conv(A)\times \{0\}$  and $\conv(A)\times \{1\}$.

It is easy to check, even by hand, that each of the simplices listed
above is unimodular when
regarded in $\A$. For example, consider the first one:
\[
\bordermatrix
{
&  \aa_0  &  \aa_4  &  \aa_1  &  \aa_2  &  \bb_2  &\cc_{+,+,+,0} \cr
&    0    &    0    & {\bf 1} &    0    &    0    &    1    \cr
&    0    &    0    &    0    & {\bf 1} &    1    &    1    \cr
&    0    &    0    &    0    &    0    &    0    & {\bf 1} \cr
&    0    & {\bf 1} &    0    &    0    &    0    &    0    \cr
&    0    &    0    &    0    &    0    & {\bf 1} &    1    \cr
& {\bf 1} &    1    &    1    &    1    &    1    &    2    \cr
}
\]
Its determinant is clearly $\pm 1$ because only the highlighted entries
produce a non-zero
summand in its expansion. The same occurs with the next 25 simplices in Table
\ref{table.trian26}.
The last two simplices in the list
are related to the previous-to-last two by the unimodular
transformation
$(x_1,x_2,x_3,x_4,x_5,x_6)\mapsto(x_1$, $x_2$, $x_3$, $x_4$, $x_6-x_5,x_6)$.
\end{proof}

Observe that only the two points $a_0$ and $b_0$ in $A_{26}$ are not
vertices. As we did with $A_{50}$ we can pull them out to become vertices,
keeping the triangulation we have constructed:

\begin{proposition}
\label{prop.26convex}
Let $A'_{26}$ be the point set obtained from $A_{26}$ by moving the
points
$a_0=(O,0)$ and  $b_0=(O,1)$ to $(O,\alpha)$ and $(O,\beta)$, 
where $\alpha\in(-\epsilon,0)$ and
$\beta\in(1,\epsilon)$ for a sufficiently small positive $\epsilon$.

Then, $A'_{26}$ is in convex position and its graph of
triangulations has at least 17 connected components.
\end{proposition}

\begin{proof}
The only difference 
with the proof of Proposition \ref{prop.50convex}
is that now 
the points $(O,\alpha)$ and $(O,\beta)$ only 
lie beyond the facets $A\times\{0\}$
and $A\times\{1\}$ for a sufficiently small perturbation.
\end{proof}

\end{document}